\documentclass[10pt]{article}
\usepackage{color}
\usepackage{amssymb, amsmath}
\usepackage{epsfig}
\usepackage{graphicx}

\newcommand{\bbc}{{\mathcal C}}

 \newcommand{\bbR}{{\mathbb R}}

\newtheorem{theorem}{Theorem}
 
 \newtheorem{proposition}{Proposition}
\newtheorem{lemma}{Lemma}%

\newcommand{\proof}{\noindent\textbf{Proof.~}}
\newcommand{\vep}{\varepsilon}
\newcommand{\qed}{\space\hfill\hspace*{\fill} $\vbox{\hrule\hbox{\vrule
height1.3ex\hskip1.3ex\vrule}\hrule}$\hss\vskip\topsep\relax}

\begin{document}

\title{On  explicit numerical schemes for the CIR process}

\author{Nikolaos Halidias \\
{\small\textsl{Department of Mathematics }}\\
{\small\textsl{University of the Aegean }}\\
{\small\textsl{Karlovassi  83200  Samos, Greece} }\\
{\small\textsl{email: nikoshalidias@hotmail.com}}}

\maketitle

\begin{abstract} In this paper we generalize an  explicit numerical
scheme for the CIR process that we have proposed before. The
advantage of the new proposed scheme is that preserves positivity
 and is well posed for a (little bit) broader set of parameters among the positivity preserving schemes. The order of convergence is at least logarithmic in general and
for a smaller set of parameters is at least $1/4$. Next we give a
different explicit numerical scheme based on exact simulation and
 we use this idea to approximate the two factor CIR model.
 Finally, we give a second explicit numerical scheme for the two
 factor CIR model based on the idea of the second section.
\end{abstract}



\section{Introduction}
Let $(\Omega, {\cal F}, \mathbb{P}, {\cal F}_t)$ be a complete
probability space with a filtration and let  a Wiener process
$(W_t)_{t \geq 0}$ defined on this space. We consider here the CIR
process, (see \cite{Cox}),
\begin{eqnarray}
x_t = x_0 + \int_0^t (kl - k x_s)ds + \sigma \int_0^t \sqrt{x_s}
dW_s,
\end{eqnarray}
where $k,l,\sigma \geq 0$. It is well known that this sde has a
unique strong solution which remain nonnegative. This stochastic
process is widely used in financial mathematics. It is well known
that one can use exact simulation methods to construct the true
solution but the drawback of such an approach is the computational
time that requires. Therefore, many researchers work on
construction of fast and efficient methods to approximate this
process. In \cite{Higham} the authors proposed a modified Euler
scheme for the approximation of the CIR process. However, this
scheme does not preserve positivity which is a desirable property
in some cases. Next, in \cite{Alfonsi}, the author proposes a
positivity preserving numerical scheme which is strongly
convergent but not for all possible parameters. In
\cite{Halidias1} we have proposed another positivity preserving
numerical scheme for the CIR process and our goal here is to
propose a generalization of this scheme in order to be well posed
for a broader class of parameters being of course positivity
preserving scheme.

 Let $0 = t_0 <
t_1 < ...<t_n = T$ and set $\Delta = \frac{T}{n}$. Consider the
following stochastic process
\begin{eqnarray}
y_t = \left( \frac{\sigma}{2(1+ka\Delta)} (W_t - W_{t_k}) +
\sqrt{y_{t_k}(1-\frac{k\Delta}{1+ka\Delta})
+\frac{\Delta}{1+ka\Delta} (kl - \frac{\sigma^2}{4(1+ka\Delta)})}
\right)^2 = (z_t)^2,
\end{eqnarray}
for $t \in (t_k,t_{k+1}]$ and a  parameter $a \in [0,1] $ where
\begin{eqnarray*}
z_t = \frac{\sigma}{2(1+ka\Delta)} (W_t - W_{t_k}) +
\sqrt{y_{t_k}(1-\frac{k\Delta}{1+ka\Delta})
+\frac{\Delta}{1+ka\Delta} (kl - \frac{\sigma^2}{4(1+ka\Delta)})},
\end{eqnarray*}
for $t \in (t_k,t_{k+1}]$.

Note that this process is well defined when $a\Delta \geq
\frac{\sigma^2-4kl}{4k^2l}$  and has the differential form, for $t
\in (t_k,t_{k+1}]$,
\begin{eqnarray}
y_t = & & y_{t_k}   + \Delta \left(kl -
\frac{\sigma^2}{4(1+ka\Delta)} - k(1-a) y_{t_k} - ka y_t \right) +
\int_{t_k}^t \frac{\sigma^2}{4(1+ka\Delta)} ds \nonumber \\ & &+
\sigma \int_{t_k}^t sgn(z_s) \sqrt{y_s} dW_s.
\end{eqnarray}
To obtain the above form we first use Ito's formula on $y_t$ and
then some simple rearrangements. This stochastic process is not
continuous in all $[0,T]$, because there are jumps at the nodes
$t_k$.

The numerical scheme that we propose here to approximate the CIR
process is the following,
\begin{eqnarray*}
y_{t_{k+1}} = \left( \frac{\sigma}{2(1+ka\Delta)} (W_{t_{k+1}} -
W_{t_k}) + \sqrt{y_{t_k}(1-\frac{k\Delta}{1+ka\Delta})
+\frac{\Delta}{1+ka\Delta} (kl - \frac{\sigma^2}{4(1+ka\Delta)})}
\right)^2,
\end{eqnarray*}
with $y_{t_0} = x_0$. Using Ito's formula one can easily see that
$y_t$ in (2) is the unique solution of the stochastic differential
equation (3). Therefore it is clear that is positivity preserving
and well defined for $a\Delta \geq \frac{\sigma^2-4kl}{4k^2l}$.
This set of parameters is (a little bit) broader than the existing
numerical schemes that preserves positivity, which usually is $4kl
\geq \sigma^2$. The main goal of future research will be the
construction of positivity preserving numerical methods that will
be well posed for all possible parameters, see for example
\cite{Andersen} for such a method without a  theoretical
convergence result.

 For  generalizations of
the semi discrete method see \cite{Halidias2}, \cite{Halidias3}.

\section{Main Results}

We will use    a  compact form of (3), for $t \in (t_k,t_{k+1}]$,
\begin{eqnarray*}
y_t =  & & x_0 + \int_0^{t} (kl  - k (1-a)y_{\hat{s}}-k a
y_{\tilde{s}})ds
\\ & & + \int_t^{t_{k+1}} \left(kl - \frac{\sigma^2}{4(1+ka\Delta)} -
k(1-a) y_{t_k} - k a y_{t} \right) ds
 + \sigma \int_0^t sgn(z_s)
\sqrt{y_s}dW_s,
\end{eqnarray*}
where
\begin{eqnarray*}
\tilde{s} = \left\{%
\begin{array}{ll}
    t_{j+1}, & \mbox{ when } s \in [t_j,t_{j+1}], \quad j=0,...,k-1\\
    t, & \mbox{ when } s \in (t_k,t]\\
\end{array}%
\right.
\end{eqnarray*}
 and     $\hat{s} = t_j$ when $s \in (t_j, t_{j+1}], j=0,...,k$.
Therefore, $y_t$ remains nonnegative as is the same as in (3).

We will remove the term $sgn(z_s)$ by changing the Brownian
motion. Set
\begin{eqnarray*}
\hat{W}_t = \int_0^t sgn(z_s) dW_s.
\end{eqnarray*}
It is easy to see that $\hat{W}$ is a continuous martingale on
${\cal F}_t$ with variation $<\hat{W},\hat{W}> = t$. Therefore,
using Levy's martingale characterization of Brownian motion (see
\cite{karatzas}, p. 157) we deduce that $\hat{W}_t$ is also a
Brownian motion. Therefore, $y_t$ satisfies the following
equation,
\begin{eqnarray*}
y_t =  & & x_0 + \int_0^{t} (kl  - k (1-a)y_{\hat{s}}-k a
y_{\tilde{s}})ds
\\ & & + \int_t^{t_{k+1}} \left(kl - \frac{\sigma^2}{4(1+ka\Delta)} -
k(1-a) y_{t_k} - k a y_{t} \right) ds
 + \sigma \int_0^t
\sqrt{y_s}d\hat{W}_s,
\end{eqnarray*}

Let now the following sde,
\begin{eqnarray}
\hat{x}_t = x_0 + \int_0^t (kl - k \hat{x}_s)ds + \sigma \int_0^t
\sqrt{\hat{x}_s} d \hat{W}_s,
\end{eqnarray}
where $\hat{W}_t$, constructed as above, is a Brownian motion
depending on $\Delta$. For each $\Delta$ the above problem has a
unique solution which has the same transition density (see
\cite{Glasserman}, p. 122), independent of $\Delta$. We will show
that $\mathbb{E}| \hat{x}_t - y_t|^2 \to 0$ as $\Delta \to 0$ and
therefore our approximation converges in the mean square sense to
a stochastic process that is equal in distribution to the unique
solution of (1). We will denote $\hat{W},\hat{x}$ again by $W,x$
for notation simplicity.

{\bf Assumption A} We suppose that $x_0 \geq 0$ a.s.
$\mathbb{E}x_0^p < A$ for some $p\geq 2$, $d = kl -
\frac{\sigma^2}{4(1+ka\Delta)} \geq 0$ and $ \Delta (1-a) \leq
\frac{1}{k}$.

\begin{lemma}[Moment bounds]
Under Assumption A  we have the moment bounds,
\begin{eqnarray*}
\mathbb{E}y_t^p + \mathbb{E} x_t^p  < C,
\end{eqnarray*}
for some $C > 0$
\end{lemma}

\proof Note that
\begin{eqnarray*}
0 \leq y_t \leq v_t = x_0 + Tkl + \sigma\int_0^t \sqrt{y_s} dW_s.
\end{eqnarray*}

Consider the stopping time $\theta_R = \inf \{ t \geq 0 :  v_t > R
\}$. Using Ito's formula on $v_{t \wedge \theta_R}^p$ we obtain,
\begin{eqnarray*}
v_{t \wedge \theta_R}^p = (x_0+Tkl)^p +
\frac{p(p-1)}{2}\sigma^2\int_0^t v_{s \wedge \theta_R}^{p-2} y_{s
\wedge \theta_R} ds + p\sigma \int_0^t v_{s \wedge \theta_R}^{p-1}
 \sqrt{y_{s \wedge \theta_R}} dW_s.
\end{eqnarray*}
Taking expectations on both  sides and noting that $y_t \leq v_t$,
we arrive at
\begin{eqnarray*}
\mathbb{E}v_{t \wedge \theta_R}^p & \leq & \mathbb{E} (x_0+Tkl)^p
+ \frac{p(p-1)}{2} \sigma^2 \int_0^t \mathbb{E} v_{s \wedge
\theta_R}^{p-1} ds \\ & \leq & \mathbb{E} (x_0+Tkl)^p +
\frac{p(p-1)}{2} \sigma^2 \int_0^t (\mathbb{E} v_{s \wedge
\theta_R}^p)^{\frac{p-1}{p}} ds
\end{eqnarray*}
Using now a Gronwall type theorem  (see \cite{Mitrinovic}, Theorem
1, p. 360), we arrive at
\begin{eqnarray}
\mathbb{E} v_{t \wedge \theta_R}^p \leq
\left([\mathbb{E}(x_0+Tkl)^p]^{\frac{p-1}{p}}+\frac{T}{2}(p-1)\sigma^2\right)^{\frac{p}{p-1}}.
\end{eqnarray}
 But $\mathbb{E} v_{t \wedge
\theta_R}^p = \mathbb{E} (v_{t \wedge \theta_R}^p \mathbb{I}_{ \{
\theta_R \geq t \} }) +
 R^p P ( \theta_R < t  )$. That means that $P( t \wedge \theta_R < t) = P ( \theta_R < t  ) \to 0$ as
 $R \to \infty$ so $t \wedge \theta_R \to t$ in probability and noting that $\theta_R$ increases as $R$ increases we have that $t \wedge \theta_R \to t$
  almost surely too, as $R \to \infty$. Going back to (4) and
 using  Fatou's lemma  we obtain,
 \begin{eqnarray*}
\mathbb{E} v_t^p \leq
\left([\mathbb{E}(x_0+Tkl)^p]^{\frac{p-1}{p}}+\frac{T(p-1)\sigma^2}{2}\right)^{\frac{p}{p-1}}
\end{eqnarray*}
We have assume in our assumptions that $\mathbb{E}x_0^p < \infty$
in order the term $\mathbb{E} (x_0+Tkl)^p$ to be well posed.

 The same holds for $x_t$ (see for example \cite{Dufresne}).
 \qed

Consider the auxiliary stochastic process, for $t \in
(t_k,t_{k+1}]$,
\begin{eqnarray}
h_t = x_0 + \int_0^t (kl  - k (1-a)y_{\hat{s}}-k a
y_{\tilde{s}})ds + \sigma \int_0^t  \sqrt{y_s}dW_s,
\end{eqnarray}
where $\tilde{s}, \hat{s}$ defined as before.

\begin{lemma}
We have the following estimates,
\begin{eqnarray*}
\mathbb{E}|h_s-y_s|^2 & \leq &  C_1 \Delta^2  \mbox{ for any } s
\in [0,T]
 \\
\mathbb{E} |h_s - y_{r}|^2 & \leq & C_2 \Delta   \mbox{ when } s
\in [t_k,t_{k+1}] \mbox{ and } r=t_k \mbox{ or }  t_{k+1}  \\
\mathbb{E}|h_s|^2 & < &  A,   \mbox{ for any } s \in [0,T].
\end{eqnarray*}
\end{lemma}

\proof Using the moment bound for $y_t$ we  easily obtain the fact
that
\begin{eqnarray*}
\mathbb{E}|h_t - y_t |^2 \leq C \Delta^2.
\end{eqnarray*}

Next, we have
\begin{eqnarray*}
\mathbb{E}|h_s - y_{t_k}|^2 \leq 2 \mathbb{E} |h_s - y_s|^2 + 2
\mathbb{E}|y_s-y_{t_k}|^2 \leq C \Delta^2 + C \Delta \leq C
\Delta.
\end{eqnarray*}

Moreover,
\begin{eqnarray*}
\mathbb{E}|h_s-y_{t_{k+1}}|^2 \leq \mathbb{E} |h_s - y_s|^2 + 2
\mathbb{E}|y_s-y_{t_{k+1}}|^2 \leq C \Delta^2 + C \Delta \leq C
\Delta.
\end{eqnarray*}
Finally, to get the moment bound for $h_t$ we just use the fact
that is close to $y_t$, i.e.
\begin{eqnarray*}
\mathbb{E}h_t^2 \leq 2\mathbb{E} |h_t-y_t|^2 +  2 \mathbb{E} y_t^2
\leq C.
\end{eqnarray*}

 \qed

\begin{theorem}
If Assumption A holds then
\begin{eqnarray*}
 \mathbb{E} |x_t - y_t|^2  \leq C \frac{1}{\sqrt{\ln n}}
\end{eqnarray*}
for any $t \in [0,T]$.
\end{theorem}

\proof

Applying Ito's formula on $|x_t-h_t|^2$ we obtain
\begin{eqnarray}
\mathbb{E} |x_t-h_t|^2 \leq  & &   2 k(1-a) \int_0^t
\mathbb{E}|x_s-h_s|
 |y_{\hat{s}} - x_s|ds + 2 ka \int_0^t \mathbb{E}  |x_s-h_s||y_{\tilde{s}}-x_s|
 ds \nonumber \\ & & + \sigma^2\int_0^t  \mathbb{E}|x_s-y_s|ds
\end{eqnarray}

Let us estimate the above quantities. It is easy to see that, for
example,
\begin{eqnarray*}
\mathbb{E} |x_s-h_s| |y_{\hat{s}} - x_s| \leq  \mathbb{E}|x_s-h_s|
(|x_s-h_s| + |h_s - y_{\hat{s}}|)
\end{eqnarray*}
Therefore, we obtain, using Cauchy-Schwarz inequality
\begin{eqnarray*}
\mathbb{E}|x_s-h_s|
 |y_{\hat{s}} - x_s| \leq \mathbb{E} |x_s - h_s|^2 +
 \sqrt{\mathbb{E} |x_s-h_s|^2} \sqrt{\mathbb{E}
 |h_s-y_{\hat{s}}|^2}, \\
 \mathbb{E}  |x_s-h_s||y_{\tilde{s}}-x_s| \leq \mathbb{E} |x_s - h_s|^2 +
 \sqrt{\mathbb{E} |x_s-h_s|^2} \sqrt{\mathbb{E}
 |h_s-y_{\tilde{s}}|^2}
 \end{eqnarray*}

Summing up we arrive at
\begin{eqnarray}
\mathbb{E}|x_t - h_t|^2 \leq C \sqrt{\Delta}+ C \int_0^t
\mathbb{E}|x_s-h_s|^2 ds + \int_0^t \mathbb{E} |x_s-h_s| ds.
\end{eqnarray}

Therefore, we have to estimate $\mathbb{E}|x_t - h_t|$. Let the
non increasing sequence $\{e_m\}_{m\in\mathbb{N}}$ with
$e_m=e^{-m(m+1)/2}$ and $e_0=1.$ We introduce the following
sequence of smooth approximations of $|x|,$ (method of Yamada and
Watanabe, \cite{Yamada})
$$
\phi_m(x)=\int_0^{|x|}dy\int_0^{y}\psi_m(u)du,
$$
where the existence of the continuous function $\psi_m(u)$ with
$0\leq \psi_m(u) \leq 2/(mu)$ and support in $(e_m,e_{m-1})$ is
justified by $\int_{e_m}^{e_{m-1}}(du/u)=m.$ The following
relations hold for $\phi_m\in\bbc^2(\bbR,\bbR)$ with
$\phi_m(0)=0,$
 $$
 |x| - e_{m-1}\leq\phi_m(x)\leq |x|, \quad |\phi_{m}^{\prime}(x)|\leq1, \quad x\in\bbR, $$
 $$
 |\phi_{m}^{\prime \prime }(x)|\leq\frac{2}{m|x|}, \,\hbox{ when }  \,e_m<|x|<e_{m-1} \,\hbox{ and }  \,  |\phi_{m}^{\prime \prime }(x)|=0 \,\hbox{ otherwise. }
 $$

Applying Ito's formula on $\phi_m(x_t-h_t)$ we obtain
\begin{eqnarray*}
\mathbb{E}\phi_m(x_t-h_t) = & & \int_0^t
\mathbb{E}\phi_m^{'}(x_s-h_s)(k(1-a)(y_{\hat{s}}-x_s) + ka(y_{\tilde{s}}-x_s))ds \\
& & + \int_0^t \frac{1}{2} \mathbb{E}\phi_m^{''}(x_s-h_s) \left(
\sigma   \sqrt{y_s} - \sigma \sqrt{x_s} \right)^2 ds.
\end{eqnarray*}

We continue by estimating
\begin{eqnarray*}
 & & \mathbb{E}\phi_m^{'}(x_s-h_s) \left(k(1-a)(y_{\hat{s}}-x_s) + k a
(y_{\tilde{s}}-x_s)\right)
\\ & \leq & k\mathbb{E}|x_s-h_s| + k(1-a)\mathbb{E}|h_s-y_{\hat{s}}| + ka
\mathbb{E} |h_s - y_{\tilde{s}}|  \\ &\leq  & k\mathbb{E}|x_s-h_s|
+ C \sqrt{\Delta}.
\end{eqnarray*}

Next,
\begin{eqnarray*}
\mathbb{E}\phi_m^{''}(x_s-h_s) \left( \sigma  \sqrt{y_s} - \sigma
\sqrt{x_s} \right)^2 \leq  \frac{4 \sigma^2}{m} + \frac{4
\sigma^2}{m} \mathbb{E}\frac{|h_s-y_s|}{e_{m}}
\end{eqnarray*}
Working as before and using Lemma 2 we get
\begin{eqnarray*}
\mathbb{E}\phi_m^{''}(x_s-h_s) \left( \sigma  \sqrt{y_s} - \sigma
\sqrt{x_s} \right)^2 \leq \frac{4 \sigma^2}{m} + C
\frac{\sqrt{\Delta}}{me_m} + C \frac{\sqrt{\Delta}}{m}.
\end{eqnarray*}

Therefore,
\begin{eqnarray*}
\mathbb{E}|x_t-h_t|  \leq e_{m-1} + \frac{4 \sigma^2 T}{m} + C
\frac{\sqrt{\Delta}}{me_m} + C \frac{\sqrt{\Delta}}{m} + k
\int_0^t \mathbb{E} |x_s-h_s|ds.
\end{eqnarray*}

Use now  Gronwall's inequality and substitute  in (9) and then
again Gronwall's inequality we arrive at
\begin{eqnarray*}
\mathbb{E}|x_t-h_t|^2 \leq C \sqrt{\Delta} + C
\frac{\sqrt{\Delta}}{me_m} + e_{m-1}.
\end{eqnarray*}
Choosing $m = \sqrt{\ln n^{\frac{1}{3}}}$ we deduce that
\begin{eqnarray*}
\mathbb{E}|x_t-h_t|^2 \leq C \frac{1}{\sqrt{\ln n}}
\end{eqnarray*}

 But
\begin{eqnarray*}
\mathbb{E}|x_t-y_t|^2 \leq 2 \mathbb{E}|x_t-h_t|^2 + 2
\mathbb{E}|h_t-y_t|^2 \leq C \frac{1}{\sqrt{\ln n}}
\end{eqnarray*}

 \qed

\section{On the polynomial rate of convergence}

We study in this section the polynomial order of convergence of
our scheme. We use a stochastic time change proposed in
\cite{Berkaoui}. For simplicity, we take $a=0$.

Our result is as follows.
\begin{proposition}
If \begin{eqnarray*} \sigma^2 \leq 2kl \mbox{ and }
\frac{1}{16}(\frac{2kl}{\sigma^2} - 1)^2 > 1
\end{eqnarray*} the following rate of convergence holds, assuming
that $x_0 \in \mathbb{R}$ and $x_0 > 0$,
\begin{eqnarray*}
\mathbb{E} | x_t - y_t |^2 \leq C  \Delta.
\end{eqnarray*}
\end{proposition}

\proof

Define the process
\begin{eqnarray*}
\gamma(t) = \int_0^t \frac{ds}{(\sqrt{x_s} + \sqrt{h_s})^2},
\end{eqnarray*}
and then the stopping time defined by
\begin{eqnarray*}
\tau_l = \inf \{ s \in [0,T]: 2 \sigma^2 \gamma(s)+ 3ks \geq l \}.
\end{eqnarray*}

 Using Ito's formula on $|x_{\tau}-h_{\tau}|^2$ with $\tau$
a stopping time, we obtain
\begin{eqnarray*}
\mathbb{E}(x_{\tau}-h_{\tau})^2 & \leq &  \int_0^{\tau}
\left(2k\mathbb{E}|x_s-h_s||y_{\hat{s}}-x_s| +\sigma^2 \mathbb{E}|
\sqrt{y_s} - \sqrt{x_s}|^2 \right) ds \\ & \leq & \int_0^{\tau}
\left(2k\mathbb{E}|x_s-h_s|^2 + 2k
\mathbb{E}|x_s-h_s||h_s-y_{\hat{s}}| +  \sigma^2
\mathbb{E}|\sqrt{y_s}-\sqrt{x_s}|^2 \right) ds \\ & \leq  &  3k
\int_0^{\tau} \mathbb{E} |x_s-h_s|^2ds +  \sigma^2 \int_0^{\tau}
\mathbb{E}|\sqrt{y_s}-\sqrt{x_s}|^2 ds +C \Delta.
\end{eqnarray*}

Now we work on
\begin{eqnarray*}
\int_0^{\tau} \mathbb{E}|\sqrt{y_s}-\sqrt{x_s}|^2 ds  \leq
\int_0^{\tau} 2\mathbb{E} |\sqrt{x_s}-\sqrt{h_s}|^2 + 2
\mathbb{E}|\sqrt{h_s}-\sqrt{y_s}|^2 ds \leq 2 \int_0^{\tau}
\mathbb{E} |\sqrt{x_s}-\sqrt{h_s}|^2 ds + C\Delta.
\end{eqnarray*}
But
\begin{eqnarray*}
\int_0^{\tau} \mathbb{E} |\sqrt{x_s}-\sqrt{h_s}|^2 ds =
\int_0^{\tau} \mathbb{E}
\frac{|x_s-h_s|^2}{(\sqrt{x_s}+\sqrt{h_s})^2} ds.
\end{eqnarray*}
Therefore,
\begin{eqnarray}
\mathbb{E}(x_{\tau} - h_{\tau})^2 \leq C\Delta +\mathbb{E}
\int_0^{\tau} |x_s-h_s|^2(3ks+2 \sigma^2 \gamma_s)^{'} ds.
\end{eqnarray}
Now, for $\tau = \tau_l$, we use the change of variables setting
$u = 3ks + 2 \sigma^2 \gamma_s$ and therefore $s = \tau_u$
obtaining,
\begin{eqnarray*}
\mathbb{E} (x_{\tau_l}-h_{\tau_l})^2 \leq C \Delta +  \int_0^l
\mathbb{E} |x_{\tau_u}-h_{\tau_u}|^2 du.
\end{eqnarray*}
Using Gronwall's inequality we obtain,
\begin{eqnarray}
\mathbb{E}|x_{\tau_l} - h_{\tau_l}|^2 \leq C e^l \Delta.
\end{eqnarray}
Going back to (9), for $\tau = t \in [0,T]$, we have under the
change of variables $u = 2 \sigma^2 \gamma_s + 3ks$,
\begin{eqnarray}
\mathbb{E}(x_t - h_t)^2 & \leq & C\Delta +\mathbb{E} \int_0^{3kT +
2 \sigma^2 \gamma_T} |x_{\tau_u}-h_{\tau_u}|^2 du \nonumber \\  &
\leq & C \Delta + \int_0^{\infty} \mathbb{E} \left( \mathbb{I}_{
\{ 3kT + 2 \sigma^2 \gamma_T \geq u \} } |x_{\tau_u} -
h_{\tau_u}|^2 \right) du.
\end{eqnarray}
Noting that
\begin{eqnarray*}
& & \int_0^{\infty} \mathbb{E} \left( \mathbb{I}_{ \{ 3kT + 2
\sigma^2 \gamma_T \geq u \} } |x_{\tau_u} - h_{\tau_u}|^2 \right)
du
\\ & \leq  & \int_0^{3kT} \mathbb{E} |x_{\tau_u} - h_{\tau_u}|^2du
+ \int_{3kT}^{\infty} \mathbb{P} (3kT + 2 \sigma^2 \gamma_T \geq
u) \mathbb{E} \left(|x_{\tau_u} - h_{\tau_u}|^2 \; | \; \{ 3kT + 2
\sigma^2
\gamma_T \geq u \} \right) du,   \\
& \leq & C \Delta + \int_{0}^{\infty} \mathbb{P} ( 2 \sigma^2
\gamma_T \geq u) \mathbb{E} |x_{\tau_u} - h_{\tau_u}|^2 du
\end{eqnarray*}
and then  we arrive using (10)
\begin{eqnarray*}
\mathbb{E}(x_t - h_t)^2 \leq C\Delta \left(1+  C \int_0^{\infty}
\mathbb{P} ( 2 \sigma^2 \gamma_T \geq u) e^{u} du \right).
\end{eqnarray*}
We will estimate now the following,
\begin{eqnarray*}
\mathbb{P}( 2 \sigma^2 \gamma_T \geq u) \leq \frac{1}{e^{mu}}
\mathbb{E}(e^{2\sigma^2 m \gamma_T} )
\end{eqnarray*}
Choose  $m = \frac{1}{16}(\frac{2kl}{\sigma^2} - 1)^2$ and use
Thm. 3.1 of \cite{Kuznetsov} to end the proof.
 \qed

In order to avoid the difficulties from the appearance of the term
$sgn(z_t)$ we have changed the Brownian motion. Below, we give a
lemma which one can use to prove strong convergence without
changing the Brownian motion and the  difference is that the order
of convergence is, at least, $\Delta^{1/4-\vep}$ for any $\vep >
0$. We prove it for the case $a=0$ for simplicity but the same
result holds for any $a \in [0,1]$.

\begin{lemma}
We have the following estimate,
\begin{eqnarray*}
\mathbb{E} y_t \left(sgn(z_t) - 1\right)^2 \leq C
\Delta^{\frac{1}{2} - \vep},
\end{eqnarray*}
for any $\vep > 0$.
\end{lemma}

\proof We begin with, when $t \in [t_k,t_{k+1}]$,
\begin{eqnarray*}
\mathbb{E} y_t (sgn(z_t) - 1)^2 & = & 4\mathbb{E}y_t \mathbb{I}_{
\{
 z_t \leq 0 \} } \leq 4 \mathbb{E} |y_t-y_{t_k}| + 4 \mathbb{E} y_{t_k} \mathbb{I}_{
\{
 z_t \leq 0 \} }   \\ & \leq & C \Delta +4 \mathbb{E} y_{t_k} \mathbb{I}_{
\{
 z_t \leq 0 \} } \mathbb{I}_{ \{ y_{t_k} \leq \Delta^{1/2-\vep} \} } +
 4 \mathbb{E} y_{t_k} \mathbb{I}_{
\{
 z_t \leq 0 \} } \mathbb{I}_{ \{ y_{t_k} > \Delta^{1/2-\vep} \} } \\
 & \leq & C \Delta^{1/2-\vep} + 4 \mathbb{E} y_{t_k} \mathbb{I}_{
\left\{
 \{z_t \leq 0 \} \cap \{  y_{t_k} > \Delta^{1/2-\vep} \} \right\} }
 \end{eqnarray*}
 We have used Lemma 2 to obtain the third inequality, estimating  the term $\mathbb{E} |y_t-y_{t_k}|$.
But
\begin{eqnarray*}
 \{ z_t \leq 0\} \cap \{ y_{t_k} > \Delta^{1/2-\vep } \}
& = & \left\{ W_t - W_{t_k} \leq -\frac{2}{\sigma} \sqrt{y_{t_k}
(1-k\Delta) + \Delta (kl - \frac{\sigma^2}{4}) } \right\} \cap \{ y_{t_k} > \Delta^{1/2-\vep } \} \\
& \subseteq & \left\{ W_t - W_{t_k} \leq -\frac{2}{\sigma}
\sqrt{1-k\Delta} \sqrt{\Delta^{1/2-\vep}} \right\} .
\end{eqnarray*}
 Since the increment $W_t-W_{t_k}$ is normally distributed
with mean zero and variance $t-t_k$ we  have that
\begin{eqnarray*}
\mathbb{P}\left( \{ z_t \leq 0 \} \cap \{ y_{t_k} >
\Delta^{1/2-\vep } \} \right) \leq \frac{\sqrt{t-t_k}}{\sqrt{2 \pi
(t-t_k)}} \int_{\frac{2 \sqrt{1-k\Delta}
\sqrt{\Delta^{1/2-\vep}}}{\sqrt{t-t_k}}}^{\infty} e^{-y^2/2} dy
\leq  \frac{C\Delta^{\vep}}{e^{C/\Delta^{\vep}}}.
\end{eqnarray*}
We have used the inequality of problem 9.22, p.112 of
\cite{karatzas} to obtain the last inequality. Now we have, using
the moment bounds for the numerical solution,
\begin{eqnarray*}
\mathbb{E} y_{t_k} \mathbb{I}_{ \left\{
 \{z_t \leq 0 \} \cap \{  y_{t_k} > \Delta^{1/2-\vep} \} \right\}
 }  \leq C \mathbb{P}\left( \{ z_t \leq 0 \} \cap \{ y_{t_k} >
\Delta^{1/2-\vep } \} \right)
\end{eqnarray*}
Noting that  $\frac{\Delta^{\vep}}{e^{1/\Delta^{\vep}}} \to 0$
faster than any power of $\Delta$  we have that
\begin{eqnarray*}
\mathbb{E} y_t (sgn(z_t) - 1)^2 \leq C \Delta^{\frac{1}{2} -
\vep}.
\end{eqnarray*}
 \qed

\section{An explicit scheme for the CIR process using exact simulation}

Consider the following equation.
\begin{eqnarray*}
x_t = x_0 + \int_0^t (kl - k x_s) ds + \int_0^t \sigma \sqrt{x_s}
dW_s.
\end{eqnarray*}
Our starting point is the exact simulation for the CIR process for
some specific parameters. If $d = \frac{4kl}{\sigma^2} \in
\mathbb{N}$ then we can simulate this process exactly (see
\cite{Glasserman}), p. 133). Indeed, the exact simulation is given
by
\begin{eqnarray*}
r(t_{i+1}) = \sum_{j=1}^{d} \left( e^{-\frac{1}{2}k\Delta}
\sqrt{\frac{r(t_i)}{d}} + \frac{\sigma}{2} \sqrt{ \frac{1}{k}
(-e^{-k\Delta})} Z_{i+1}^{(j)} \right)^2,
\end{eqnarray*}
where $(Z_i^{(1)},...,Z_i^{(d)})$ are standard normal $d$-vectors,
independent for different values of $i$. Therefore, the idea (see
\cite{Halidias1}) is to split a part of the drift term and the
remaining drift coefficient will be such that we can simulate it
exactly. Then, we will study the error produced by this splitting.
First we assume that $d > 1$ and we will propose an explicit
numerical scheme that preserves positivity and converges in the
mean square sense with, at least,  logarithmic order. For the case
$2kl > 5 \sigma^2$ we will show that this solution converges in
the mean square sense with $1/2$ order of convergence.

\subsection{The general case $d > 1$}
We will use the main idea of \cite{Halidias1} and propose the
following semi discrete numerical scheme,
\begin{eqnarray*}
y_t = y_{t_k} + \Delta k_1 l - y_{t_k} \Delta k_1 + \int_{t_k}^t
(k_2 l - k_2 y_s)ds + \sigma \int_{t_k}^t \sqrt{y_s} dW_s,
\end{eqnarray*}
where $k = k_1 + k_2$ and $\frac{4k_2l}{\sigma^2} =
[\frac{4kl}{\sigma^2}]$ and by $[\cdot]$ we denote the integer
part. The above sde has a unique strong solution which can be
simulated exactly and is well posed when $\Delta < \frac{1}{k_1}$.
A compact form of the numerical scheme is,
\begin{eqnarray*}
y_t = x_0 + \int_0^t(kl - k_2 y_s - k_1 y_{\hat{s}} ) ds +
\int_t^{t_{k+1}} (k_1 l - k_1 y_{\hat{s}})ds + \sigma \int_0^t
\sqrt{y_s} dW_s, \quad t \in (t_k,t_{k+1}].
\end{eqnarray*}

\begin{lemma}[Moment bounds]
Under Assumption A  we have the moment bounds,
\begin{eqnarray*}
\mathbb{E}y_t^p + \mathbb{E} x_t^p  < C,
\end{eqnarray*}
for some $C > 0$
\end{lemma}

\proof Note that
\begin{eqnarray*}
0 \leq y_t \leq v_t = x_0 + Tkl + \sigma\int_0^t \sqrt{y_s} dW_s.
\end{eqnarray*}

Consider the stopping time $\theta_R = \inf \{ t \geq 0 :  v_t > R
\}$. Using Ito's formula on $v_{t \wedge \theta_R}^p$ we obtain,
\begin{eqnarray*}
v_{t \wedge \theta_R}^p = (x_0+Tkl)^p +
\frac{p(p-1)}{2}\sigma^2\int_0^t v_{s \wedge \theta_R}^{p-2} y_{s
\wedge \theta_R} ds + p\sigma \int_0^t v_{s \wedge \theta_R}^{p-1}
 \sqrt{y_{s \wedge \theta_R}} dW_s.
\end{eqnarray*}
Taking expectations on both  sides and noting that $y_t \leq v_t$,
we arrive at
\begin{eqnarray*}
\mathbb{E}v_{t \wedge \theta_R}^p & \leq & \mathbb{E} (x_0+Tkl)^p
+ \frac{p(p-1)}{2} \sigma^2 \int_0^t \mathbb{E} v_{s \wedge
\theta_R}^{p-1} ds \\ & \leq & \mathbb{E} (x_0+Tkl)^p +
\frac{p(p-1)}{2} \sigma^2 \int_0^t (\mathbb{E} v_{s \wedge
\theta_R}^p)^{\frac{p-1}{p}} ds
\end{eqnarray*}
Using now a Gronwall type theorem  (see \cite{Mitrinovic}, Theorem
1, p. 360), we arrive at
\begin{eqnarray}
\mathbb{E} v_{t \wedge \theta_R}^p \leq
\left([\mathbb{E}(x_0+Tkl)^p]^{\frac{p-1}{p}}+\frac{T}{2}(p-1)\sigma^2\right)^{\frac{p}{p-1}}.
\end{eqnarray}
 But $\mathbb{E} v_{t \wedge
\theta_R}^p = \mathbb{E} (v_{t \wedge \theta_R}^p \mathbb{I}_{ \{
\theta_R \geq t \} }) +
 R^p P ( \theta_R < t  )$. That means that $P( t \wedge \theta_R < t) = P ( \theta_R < t  ) \to 0$ as
 $R \to \infty$ so $t \wedge \theta_R \to t$ in probability and noting that $\theta_R$ increases as $R$ increases we have that $t \wedge \theta_R \to t$
  almost surely too, as $R \to \infty$. Going back to (4) and
 using  Fatou's lemma  we obtain,
 \begin{eqnarray*}
\mathbb{E} v_t^p \leq
\left([\mathbb{E}(x_0+Tkl)^p]^{\frac{p-1}{p}}+\frac{T(p-1)\sigma^2}{2}\right)^{\frac{p}{p-1}}
\end{eqnarray*}
We have assume in our assumptions that $\mathbb{E}x_0^p < \infty$
in order the term $\mathbb{E} (x_0+Tkl)^p$ to be well posed.

 The same holds for $x_t$.
 \qed

Consider now the following auxiliary stochastic process,
\begin{eqnarray*}
h_t = x_0 +\int_0^t(kl - k_2 y_s - k_1 y_{\hat{s}} ) ds  + \sigma
\int_0^t \sqrt{y_s} dW_s, \quad t \in (t_k,t_{k+1}].
\end{eqnarray*}

\begin{lemma}
We have the following estimates,
\begin{eqnarray*}
\mathbb{E}|h_s-y_s|^2 & \leq &  C_1 \Delta^2  \mbox{ for any } s
\in [0,T]
 \\
\mathbb{E} |h_s - y_{t_k}|^2 & \leq & C_2 \Delta   \mbox{ when } s
\in [t_k,t_{k+1}]   \\
\mathbb{E}|h_s|^2 & < &  A,   \mbox{ for any } s \in [0,T].
\end{eqnarray*}
\end{lemma}
\proof Noting that
\begin{eqnarray*}
h_t - y_t = \int_t^{t_{k+1}} (k_1l - k_1 y_{\hat{s}})ds
\end{eqnarray*}
we can easily take the results. \qed

\begin{theorem}
If Assumption A holds then
\begin{eqnarray*}
 \mathbb{E} |x_t - y_t|^2  \leq C \frac{1}{\sqrt{\ln n}}
\end{eqnarray*}
for any $t \in [0,T]$.
\end{theorem}

\proof

Applying Ito's formula on $|x_t-h_t|^2$ we obtain
\begin{eqnarray}
\mathbb{E} |x_t-h_t|^2 \leq   \int_0^t \mathbb{E}\left( 2k_1
|h_s-x_s||x_s - y_{\hat{s}}| + 2k_2 |h_s-x_s| |x_s-y_s| + \sigma^2
|x_s-y_s| \right) ds
\end{eqnarray}

Let us estimate the above quantities. It is easy to see that,
using Young inequality,
\begin{eqnarray*}
\mathbb{E} |x_s-h_s| |y_{\hat{s}} - x_s| + \mathbb{E} |x_s-h_s|
|y_{s} - x_s|\leq C \mathbb{E}|x_s-h_s|^2 + C \sqrt{\Delta}
\end{eqnarray*}

Summing up we arrive at
\begin{eqnarray}
\mathbb{E}|x_t - h_t|^2 \leq C \sqrt{\Delta}+ C \int_0^t
\mathbb{E}|x_s-h_s|^2 ds + \sigma^2\int_0^t \mathbb{E} |x_s-h_s|
ds.
\end{eqnarray}

Therefore, we have to estimate $\mathbb{E}|x_t - h_t|$. Let the
non increasing sequence $\{e_m\}_{m\in\mathbb{N}}$ with
$e_m=e^{-m(m+1)/2}$ and $e_0=1.$ We introduce the following
sequence of smooth approximations of $|x|,$ (method of Yamada and
Watanabe, \cite{Yamada})
$$
\phi_m(x)=\int_0^{|x|}dy\int_0^{y}\psi_m(u)du,
$$
where the existence of the continuous function $\psi_m(u)$ with
$0\leq \psi_m(u) \leq 2/(mu)$ and support in $(e_m,e_{m-1})$ is
justified by $\int_{e_m}^{e_{m-1}}(du/u)=m.$ The following
relations hold for $\phi_m\in\bbc^2(\bbR,\bbR)$ with
$\phi_m(0)=0,$
 $$
 |x| - e_{m-1}\leq\phi_m(x)\leq |x|, \quad |\phi_{m}^{\prime}(x)|\leq1, \quad x\in\bbR, $$
 $$
 |\phi_{m}^{\prime \prime }(x)|\leq\frac{2}{m|x|}, \,\hbox{ when }  \,e_m<|x|<e_{m-1} \,\hbox{ and }  \,  |\phi_{m}^{\prime \prime }(x)|=0 \,\hbox{ otherwise. }
 $$

Applying Ito's formula on $\phi_m(x_t-h_t)$ we obtain
\begin{eqnarray*}
\mathbb{E}\phi_m(x_t-h_t) \leq & & \int_0^t
\mathbb{E}\phi_m^{'}(x_s-h_s)(k_1(y_{\hat{s}}-x_s) + k_2(y_{s}-x_s))ds \\
& & + \int_0^t \frac{\sigma^2}{2} \mathbb{E}\phi_m^{''}(x_s-h_s)
|x_s-y_s|  ds.
\end{eqnarray*}

We continue by estimating
\begin{eqnarray*}
  \mathbb{E}\phi_m^{'}(x_s-h_s)(k_1(y_{\hat{s}}-x_s) + k_2(y_{s}-x_s))
   \leq  C\mathbb{E}|x_s-h_s|
+ C \sqrt{\Delta}.
\end{eqnarray*}

Next,
\begin{eqnarray*}
\mathbb{E}\phi_m^{''}(x_s-h_s) \left(  \sqrt{y_s} - \sigma
\sqrt{x_s} \right)^2 \leq  \frac{4 \sigma^2}{m} + \frac{4
\sigma^2}{m} \mathbb{E}\frac{|h_s-y_s|}{e_{m}}
 \leq \frac{4 \sigma^2}{m} + \frac{C}{m} \frac{\sqrt{\Delta}}{e_{m}}
\end{eqnarray*}

Therefore,
\begin{eqnarray*}
\mathbb{E}|x_t-h_t|  \leq e_{m-1} + \frac{4 \sigma^2}{m} + C
\frac{\sqrt{\Delta}}{me_m} + k \int_0^t \mathbb{E} |x_s-h_s|ds.
\end{eqnarray*}

Use now  Gronwall's inequality and substitute  in (3) and then
again Gronwall's inequality we arrive at
\begin{eqnarray*}
\mathbb{E}|x_t-h_t|^2 \leq C \sqrt{\Delta} + C
\frac{\sqrt{\Delta}}{me_m} + e_{m-1}.
\end{eqnarray*}
Choosing $m = \sqrt{\ln n^{\frac{1}{3}}}$ we deduce that
\begin{eqnarray*}
\mathbb{E}|x_t-h_t|^2 \leq C \frac{1}{\sqrt{\ln n}}
\end{eqnarray*}

 But
\begin{eqnarray*}
\mathbb{E}|x_t-y_t|^2 \leq 2 \mathbb{E}|x_t-h_t|^2 + 2
\mathbb{E}|h_t-y_t|^2 \leq C \frac{1}{\sqrt{\ln n}}
\end{eqnarray*}
 \qed

\subsection{The case $2kl > 5 \sigma^2$}

Here, we choose again $k_1,k_2$ such that $k = k_1+k_2$ and $d =
\frac{4k_2l}{\sigma^2} = [ \frac{4kl}{\sigma^2} ]$. Our result is
as follows.
\begin{proposition}
If \begin{eqnarray*} \sigma^2 \leq 2kl \mbox{ and }
\frac{1}{16}(\frac{2kl}{\sigma^2} - 1)^2 > 1
\end{eqnarray*} the following rate of convergence holds, assuming
that $x_0 \in \mathbb{R}_+$,
\begin{eqnarray*}
\mathbb{E} | x_t - y_t |^2 \leq C  \Delta.
\end{eqnarray*}
\end{proposition}

\proof

Define the process
\begin{eqnarray*}
\gamma(t) = \int_0^t \frac{ds}{(\sqrt{x_s} + \sqrt{h_s})^2},
\end{eqnarray*}
and then the stopping time defined by
\begin{eqnarray*}
\tau_l = \inf \{ s \in [0,T]: 2\sigma^2\gamma(s)+ 3ks \geq l \}.
\end{eqnarray*}

 Using Ito's formula on $|x_{\tau}-h_{\tau}|^2$ with $\tau$
a stopping time, we obtain
\begin{eqnarray*}
\mathbb{E}(x_{\tau}-h_{\tau})^2 & \leq &  \int_0^{\tau}
\left(2k_1\mathbb{E}|x_s-h_s||y_{\hat{s}}-x_s| +2k_2 \mathbb{E}
|h_s-x_s| |x_s-y_s| +\sigma^2 \mathbb{E}| \sqrt{y_s} -
\sqrt{x_s}|^2 \right) ds
\\ & \leq & C \Delta + \int_0^{\tau} \left(3k\mathbb{E}|x_s-h_s|^2  +  \sigma^2 \mathbb{E}|\sqrt{y_s}-\sqrt{x_s}|^2 \right) ds
\end{eqnarray*}

Now we work on
\begin{eqnarray*}
\int_0^{\tau} \mathbb{E}|\sqrt{y_s}-\sqrt{x_s}|^2 ds  \leq
\int_0^{\tau} 2\mathbb{E} |\sqrt{x_s}-\sqrt{h_s}|^2 + 2
\mathbb{E}|\sqrt{h_s}-\sqrt{y_s}|^2 ds \leq 2 \int_0^{\tau}
\mathbb{E} |\sqrt{x_s}-\sqrt{h_s}|^2 ds + C\Delta.
\end{eqnarray*}
But
\begin{eqnarray*}
\int_0^{\tau} \mathbb{E} |\sqrt{x_s}-\sqrt{h_s}|^2 ds =
\int_0^{\tau} \mathbb{E}
\frac{|x_s-h_s|^2}{(\sqrt{x_s}+\sqrt{h_s})^2} ds.
\end{eqnarray*}
Therefore,
\begin{eqnarray}
\mathbb{E}(x_{\tau} - h_{\tau})^2 \leq C\Delta +\mathbb{E}
\int_0^{\tau} |x_s-h_s|^2(3ks+2 \sigma^2 \gamma_s)^{'} ds.
\end{eqnarray}
Now, for $\tau = \tau_l$, we use the change of variables setting
$u = 3ks + 2\sigma^2\gamma_s$ and therefore $s = \tau_u$
obtaining,
\begin{eqnarray*}
\mathbb{E} (x_{\tau_l}-h_{\tau_l})^2 \leq C \Delta +  \int_0^l
\mathbb{E} |x_{\tau_u}-h_{\tau_u}|^2 du.
\end{eqnarray*}
Using Gronwall's inequality we obtain,
\begin{eqnarray}
\mathbb{E}|x_{\tau_l} - h_{\tau_l}|^2 \leq C e^l \Delta.
\end{eqnarray}
Going back to (4), for $\tau = t \in [0,T]$, we have under the
change of variables $u = 2 \sigma^2\gamma_s + 3ks$,
\begin{eqnarray}
\mathbb{E}(x_t - h_t)^2 & \leq & C\Delta +\mathbb{E} \int_0^{3kT +
2 \sigma^2 \gamma_T} |x_{\tau_u}-h_{\tau_u}|^2 du \nonumber \\  &
\leq & C \Delta + \int_0^{\infty} \mathbb{E} \left( \mathbb{I}_{
\{ 3kT + 2 \sigma^2 \gamma_T \geq u \} } |x_{\tau_u} -
h_{\tau_u}|^2 \right) du.
\end{eqnarray}
Noting that
\begin{eqnarray*}
& & \int_0^{\infty} \mathbb{E} \left( \mathbb{I}_{ \{ 3kT + 2
\sigma^2 \gamma_T \geq u \} } |x_{\tau_u} - h_{\tau_u}|^2 \right)
du
\\ & \leq  & \int_0^{3kT} \mathbb{E} |x_{\tau_u} - h_{\tau_u}|^2du
+ \int_{3kT}^{\infty} \mathbb{P} (3kT + 2 \sigma^2 \gamma_T \geq
u) \mathbb{E} \left(|x_{\tau_u} - h_{\tau_u}|^2 \; | \; \{ 3kT + 2
\sigma^2 \gamma_T \geq u \} \right) du,   \\
& \leq & C \Delta + \int_{0}^{\infty} \mathbb{P} ( 2 \sigma^2
\gamma_T \geq u) \mathbb{E} |x_{\tau_u} - h_{\tau_u}|^2 du
\end{eqnarray*}
and then  we arrive using (5)
\begin{eqnarray*}
\mathbb{E}(x_t - h_t)^2 \leq C\Delta \left(1+  C \int_0^{\infty}
\mathbb{P} ( 2 \sigma^2 \gamma_T \geq u) e^{u} du \right).
\end{eqnarray*}
We will estimate now the following,
\begin{eqnarray*}
\mathbb{P}( 2 \sigma^2 \gamma_T \geq u) \leq \frac{1}{e^{mu}}
\mathbb{E}(e^{2m \sigma^2 \gamma_T} )
\end{eqnarray*}
Choose  $m = \frac{1}{16}(\frac{2kl}{\sigma^2} - 1)^2$ and use
Thm. 3.1 of \cite{Kuznetsov} to end the proof.
 \qed

\section{An explicit scheme for the two factor CIR model based on exact simulation}

Let $(\Omega, {\cal F}, \mathbb{P}, {\cal F}_t)$ be a complete
probability space with a filtration and let two independent Wiener
processes  $(W^{1,2}_t)_{t \geq 0}$ defined on this space. Here we
consider the following two factor CIR model (see \cite{Shreve}, p.
420),
\begin{eqnarray*}
x_1(t) & = & x_1(0) + \int_0^t (k - \lambda_{11} x_1(s) +
\lambda_{12}x_2(s))ds +
\int_0^t \sigma_1 \sqrt{x_1(s)} dW^1_s,\\
x_2(t) & = & x_2(0) + \int_0^t (l - \lambda_{21} x_2(s) +
\lambda_{22} x_1(s))ds + \int_0^t \sigma_2 \sqrt{x_2(s)} dW^2_s
\end{eqnarray*}
This kind of model is widely used in financial mathematics. If one
wants to calculate complicate expressions of the solution of the
above system maybe the only way is to approximate it numerically.
In this case, the numerical scheme should be positivity preserving
and the usual Euler scheme does not have this property. For more
details about the use of this model in financial mathematics one
can see for example \cite{Shreve}.

In the following two sections we will propose two different,
explicit and positivity preserving numerical schemes.

Our starting point is the exact simulation for the CIR process for
some specific parameters. Consider the CIR process,  and let $0 =
t_0 < t_1 < ...<t_n = T$, setting $\Delta = \frac{T}{n}$,
\begin{eqnarray*}
x_t = x_0 + \int_0^t (kl - k x_s)ds + \sigma \int_0^t \sqrt{x_s}
dW_s.
\end{eqnarray*}
If $d = \frac{4kl}{\sigma^2} \in \mathbb{N}$ then we can simulate
this process exactly (see \cite{Glasserman}), p. 133). Indeed, the
exact simulation is given by
\begin{eqnarray*}
r(t_{i+1}) = \sum_{j=1}^{d} \left( e^{-\frac{1}{2}k\Delta}
\sqrt{\frac{r(t_i)}{d}} + \frac{\sigma}{2} \sqrt{ \frac{1}{k}
(-e^{-k\Delta})} Z_{i+1}^{(j)} \right)^2,
\end{eqnarray*}
where $(Z_i^{(1)},...,Z_i^{(d)})$ are standard normal $d$-vectors,
independent for different values of $i$. Therefore, the idea (see
\cite{Halidias1}) is to split a part of the drift term and the
remaining drift coefficient will be such that we can simulate it
exactly. Then, we will study the error produced by this splitting.
For the two factor CIR model there is one more difficulty. In each
equation there exists an unknown stochastic process which appears
on the other. In this situation we will use the main idea of
\cite{Halidias2} and discretize every part of the first stochastic
differential equation that contains the unknown stochastic process
which contained in the second equation and vice versa. In this way
we arrive to two stochastic differential equations that contains
only one unknown stochastic process. For another positivity
preserving numerical scheme for one factor CIR model see
\cite{Alfonsi}.

We propose  the following decomposition,
\begin{eqnarray*}
y_1(t)  =  y_1(t_k)  & + &   \Delta \lambda_{12} y_2(t_k)  +
\Delta k_1
 +\int_{t_k}^t (k_2 - \lambda_{11} y_1(s))ds \\  & + &  \sigma_1 \int_{t_k}^t
\sqrt{y_1(s)} dW^1_s, \quad t \in (t_k,t_{k+1}]
\end{eqnarray*}
\begin{eqnarray*}
 y_2(t)  =
y_2(t_k) & + &  \Delta \lambda_{22} y_1(t_k)  +  \Delta l_1 +
\int_{t_k}^t (l_2 - \lambda_{21} y_2(s))ds \\  &  + &   \sigma_2
\int_{t_k}^t \sqrt{y_2(s)} dW^2_s \quad t \in (t_k,t_{k+1}]
\end{eqnarray*}
where $\frac{4k_2}{\sigma_1^2} = [\frac{4k}{\sigma_1^2}]$,
$\frac{4l_2}{\sigma_2^2} = [\frac{4l}{\sigma_2^2}]$ and by $[
\cdot ]$ we denote the integer part. We see that the above sdes
are not really a system and in each equation only one unknown
stochastic process appears. Therefore, in each step, we can
simulate exactly the stochastic process $y_1,y_2$.

 Let us write in a more compact form our numerical scheme, for $t
\in (t_k,t_{k+1}]$,
\begin{eqnarray*}
y_1(t)  =  x_1(0) & + & \int_0^t \left(k - \lambda_{11} y_1(s) +
\lambda_{12} y_2(\hat{s}) \right) ds +(t_{k+1}-t)(k_1+\lambda_{12}
y_2(t_k))\\ & + & \sigma_1 \int_0^t \sqrt{y_1(s)}
dW_s^1,\end{eqnarray*}
\begin{eqnarray*} y_2(t)  =   x_2(0) & + & \int_0^t \left(l -
\lambda_{21} y_2(s) + \lambda_{22} y_1(\hat{s}) \right) ds
+(t_{k+1}-t)(l_1+\lambda_{22} y_1(t_k)) \\ & + & \sigma_2 \int_0^t
\sqrt{y_2(s)} dW_s^2,
\end{eqnarray*}
where  $\hat{s} = t_k$ when $ s \in [t_k,t_{k+1}]$. Our first
result is to obtain the moment bounds for the true and the
approximate solution.

{\bf Assumption A} Assume that $x_1(0),x_2(0) \in \mathbb{R}_+$
and that $d_1 = \frac{4k}{\sigma_1^2} > 1$, $d_2 =
\frac{4l}{\sigma^2_2} > 1$.

Below we will give the moment bounds for the true and the
approximate solution. However, for the approximate solution it
seems that we need to bound it uniformly  as we did, for example
in \cite{Halidias3}.
\begin{lemma}
Under Assumption A we have
\begin{eqnarray*}
\mathbb{E}(\sup_{0 \leq t \leq T} ( y_1(t)^2 + y_2(t)^2) )< C,
\quad \mathbb{E} x_1(t)^2 + x_2(t)^2 < C.
\end{eqnarray*}
\end{lemma}

\proof We easily see that
\begin{eqnarray*}
0 \leq y_1(t) \leq v_1(t) = x_1(0) & + & Tk + \Delta \lambda_{12}
x_2(0)+\int_0^t \lambda_{12} (y_2(\hat{s})+y_2(t_k)) ds \\ &  + &
\sigma_1 \int_0^t \sqrt{y_1(s)} dW^1_s,\end{eqnarray*}
\begin{eqnarray*} 0 \leq y_2(t) \leq v_2(t) = x_2(0) & + & Tl + \Delta
\lambda_{22} x_1(0)+\int_0^t \lambda_{22} (y_1(\hat{s})+y_1(t_k))
ds \\ & + & \sigma_2\int_0^t \sqrt{y_2(s)} dW^2_s,
\end{eqnarray*}
We have used that $(t_{k+1}-t) \lambda_{12} y_2(t_k) \leq \Delta
\lambda_{12} x_2(0)$ when $t_k = t_0$ and $(t_{k+1}-t)
\lambda_{12} y_2(t_k) \leq \int_0^t \lambda_{12} y_2(t_k) ds$ when
$t_k = t_1,t_2,...$ and therefore $t > \Delta$. Thus,
\begin{eqnarray*}
(t_{k+1}-t) \lambda_{12} y_2(t_k) \leq \Delta \lambda_{12} x_2(0)
+\int_0^t \lambda_{12} y_2(t_k) ds,
\end{eqnarray*}
for any $t \in [0,T]$.
 The same holds for the $y_1(t_k)$.

Consider the stopping time $\tau = \inf \{ t \in [0,T] : y_1(t)
> R \mbox{ or } y_2(t) > R \}$. Then, we can write,
\begin{eqnarray*}
v_1^2(t \wedge \tau) \leq C + C \int_0^{t \wedge \tau} v^2_2(
\hat{s} \wedge \tau) +v_2^2(t_k \wedge \tau) ds + C |
\int_0^{t \wedge \tau} \sqrt{y_1(s \wedge \tau)}dW^1_s|^2, \\
v_2^2(t \wedge \tau) \leq C + C \int_0^{t \wedge \tau} v^2_1(
\hat{s} \wedge \tau) +v_1^2(t_k \wedge \tau) ds + C | \int_0^{t
\wedge \tau} \sqrt{y_2(s \wedge \tau)}dW^2_s|^2,
\end{eqnarray*}
and therefore
\begin{eqnarray*}
\sup_{0 \leq t \leq r} (v_1^2(t \wedge \tau) +v_2^2(t \wedge
\tau)) \leq  & & C  +   C \int_0^{r} (v_1^2(\hat{s} \wedge \tau) +
v_2^2(\hat{s} \wedge \tau)+v_2^2(t_k \wedge \tau)+v_1^2(t_k \wedge
\tau)) ds \\ &   + & C \left( \sup_{0 \leq t \leq r}|\int_0^{t
\wedge \tau} \sqrt{y_1(s)}dW^1_s|^2 + \sup_{0 \leq t \leq r}|
\int_0^{t \wedge \tau} \sqrt{y_2(s)}dW^2_s|^2 \right).
\end{eqnarray*}
Taking expectations and using Doob's martingale inequality we
arrive at
\begin{eqnarray*}
 \mathbb{E} (\sup_{0 \leq t \leq r} (v_1^2(t) +v_2^2(t)) \leq & &  C + C \int_0^{r} \mathbb{E}( (v_1^2(\hat{s} \wedge
\tau) + v_2^2(\hat{s} \wedge \tau) +v_2^2(t_k \wedge
\tau)+v_1^2(t_k \wedge \tau)) ds  \\ & + & C \int_0^{r}
\mathbb{E}(v_1(s \wedge \tau)+v_2(s \wedge \tau))ds
\end{eqnarray*} \begin{eqnarray*} \leq C + C \int_0^{r} \left( \mathbb{E} \sup_{ 0
\leq \beta \leq s} (v_1^2(\beta \wedge \tau) + v_2^2(\beta \wedge
\tau))  + \sqrt{\mathbb{E}(\sup_{0 \leq \beta \leq s} (v_1^2(\beta
\wedge \tau) + v_2^2(\beta \wedge \tau))) } \right) ds.
\end{eqnarray*}

Setting now $u(r) = \mathbb{E} (\sup_{0 \leq t \leq r} (v_1^2(t
\wedge \tau) +v_2^2(t \wedge \tau))$ and using a generalized
Gronwall inequality (see \cite{Mitrinovic}, Theorem 1, p. 360) we
deduce that
\begin{eqnarray*}
u(r) \leq C, \quad r \in [0,T]
\end{eqnarray*}
with $C$ independent of $R$. Taking the limit as $R \to \infty$
and using Fatou's lemma we take our result.
 The same holds for $x_1,x_2$.

 \qed
We will use later the auxiliary stochastic processes,
\begin{eqnarray*}
h_1(t) & = & x_1(0) + \int_0^t \left(k - \lambda_{11} y_1(s) +
\lambda_{12} y_2(\hat{s}) \right) ds  + \sigma_1
\int_0^t \sqrt{y_1(s)} dW_s^1, \\
h_2(t) & = &  x_2(0) + \int_0^t \left(l - \lambda_{21} y_2(s) +
\lambda_{22} y_1(\hat{s}) \right) ds  + \sigma_2 \int_0^t
\sqrt{y_2(s)} dW_s^2
\end{eqnarray*}
We shall show below that $h_{1,2}(t)$ and $y_{1,2}(t)$ remain
close.
\begin{lemma} Under Assumption A we have, for all $t \in [0,T]$,
\begin{eqnarray*}
\mathbb{E}|h_{1,2}(t) - y_{1,2}(t)|^2   & \leq & C \Delta^2 \\
\mathbb{E} |h_{1,2}(t) - y_{1,2}(t_k)|^2 & \leq & C \Delta \mbox{ when } t \in [t_k,t_{k+1}]\\
\mathbb{E}h_{1,2}^2(t) & \leq & C.
\end{eqnarray*}
\end{lemma}

\proof It is easy to see that
\begin{eqnarray*}
\mathbb{E}|y_{1,2}(t) - y_{1,2}(\hat{t})|^2 \leq C \Delta.
\end{eqnarray*}
Moreover, noting that
\begin{eqnarray*}
\mathbb{E}|y_{1,2}(t) - h_{1,2}(t)| \leq  C \Delta^2,
\end{eqnarray*}
we obtain the other results. \qed

\subsection{The general case $d_1 \geq  1$, $d_2 \geq  1$} In this
section we assume that $d_1 > 1$ and $d_2 > 1$
 and we will prove that the rate of convergence is at least
 logarithmic. If $d_1=1$ for example we can simulate $x_1$
 exactly therefore we work on the case where $d_1 > 1$ and $d_2 >
 1$.
\begin{theorem}
If Assumption A holds then
\begin{eqnarray*}
 \mathbb{E} (|x_1(t) - y_1(t)|^2+|x_2(t)-y_2(t)|^2)  \leq C \frac{1}{\sqrt{\ln n}}
\end{eqnarray*}
for any $t \in [0,T]$.
\end{theorem}

\proof

Applying Ito's formula on $|x_{1}(t)-h_{1}(t)|^2$  we obtain
\begin{eqnarray}
& & \mathbb{E} |x_1(t)-h_1(t)|^2 \nonumber \\ & \leq &  \int_0^t
\mathbb{E}\big( 2 \lambda_{11}|x_1(s)-h_1(s)| |y_1(s)-x_1(s)|
\nonumber  + 2\lambda_{12}|x_1(s)-h_1(s)||x_2(s)-y_2(\hat{s})| \\
& & + \sigma_1^2 |x_1(s)-y_1(s)| \big) ds
\end{eqnarray}
Using Young inequality, we deduce
\begin{eqnarray*}
& & \mathbb{E} |x_1(s)-h_1(s)| |y_1(s) - x_1(s)| + \mathbb{E}
|x_1(s)-h_1(s)| |y_2(\hat{s}) - x_2(s)| \\ & \leq & C \left(
\mathbb{E}|x_1(s)-h_1(s)|^2 + \mathbb{E}|x_2(s)-h_2(s)|^2
\right)+C \Delta
\end{eqnarray*}

Summing up we arrive at
\begin{eqnarray}
& & \mathbb{E}|x_1(t) - h_1(t)|^2 \nonumber \\ & \leq & C
\sqrt{\Delta}+ C \int_0^t
\mathbb{E}(|x_1(s)-h_1(s)|^2+|x_2(s)-h_2(s)|^2 ) ds \nonumber \\ &
& + \sigma^2_1\int_0^t \mathbb{E} |x_1(s)-h_1(s)| ds.
\end{eqnarray}

Setting $v^2(t) = |x_1(s)-h_1(s)|^2+|x_2(s)-h_2(s)|^2$, using
Ito's formula as before on $|x_2(t)-h_2(t)|^2$ and adding the
results we arrive at
\begin{eqnarray*}
\mathbb{E}v^2(t) \leq C \sqrt{\Delta}+ C \int_0^t \mathbb{E}
v^2(s)ds + (\sigma_1^2 + \sigma_2^2) \int_0^t
\mathbb{E}(|x_1(s)-h_1(s)| + |x_2(s)-h_2(s)|)ds.
\end{eqnarray*}

Therefore, we have to estimate $\mathbb{E}|x_1(t) - h_1(t)|$ and
$\mathbb{E}|x_2(t) - h_2(t)|$. Let the non increasing sequence
$\{e_m\}_{m\in\mathbb{N}}$ with $e_m=e^{-m(m+1)/2}$ and $e_0=1.$
We introduce the following sequence of smooth approximations of
$|x|,$ (method of Yamada and Watanabe, \cite{Yamada})
$$
\phi_m(x)=\int_0^{|x|}dy\int_0^{y}\psi_m(u)du,
$$
where the existence of the continuous function $\psi_m(u)$ with
$0\leq \psi_m(u) \leq 2/(mu)$ and support in $(e_m,e_{m-1})$ is
justified by $\int_{e_m}^{e_{m-1}}(du/u)=m.$ The following
relations hold for $\phi_m\in\bbc^2(\bbR,\bbR)$ with
$\phi_m(0)=0,$
 $$
 |x| - e_{m-1}\leq\phi_m(x)\leq |x|, \quad |\phi_{m}^{\prime}(x)|\leq1, \quad x\in\bbR, $$
 $$
 |\phi_{m}^{\prime \prime }(x)|\leq\frac{2}{m|x|}, \,\hbox{ when }  \,e_m<|x|<e_{m-1} \,\hbox{ and }  \,  |\phi_{m}^{\prime \prime }(x)|=0 \,\hbox{ otherwise. }
 $$

Applying Ito's formula on $\phi_m(x_1(t)-h_1(t))$ we obtain
\begin{eqnarray*}
\mathbb{E}\phi_m(x_1(t)-h_1(t)) \leq & & \int_0^t
\mathbb{E}|\phi_m^{'}(x_1(s)-h_1(s))|(\lambda_{11}|y_1(s)-x_1(s)|+\lambda_{12}|x_2(s)-y_2(\hat{s})|)ds \\
& & + \int_0^t \frac{\sigma_1^2}{2}
\mathbb{E}|\phi_m^{''}(x_1(s)-h_1(s))| |x_1(s)-y_1(s)|  ds.
\end{eqnarray*}
We continue by estimating,
\begin{eqnarray*}
 & &  \mathbb{E}|\phi_m^{'}(x_1(s)-h_1(s))|(\lambda_{11}|y_1(s)-x_1(s)|+\lambda_{12}|x_2(s)-y_2(\hat{s})|)
 \\ & \leq & C \mathbb{E}(|x_1(s)-h_1(s)|+|x_2(s)-h_2(s)|) + C
  \sqrt{\Delta},
\end{eqnarray*}
and
\begin{eqnarray*}
\mathbb{E}|\phi_m^{''}(x_1(s)-h_1(s))| |x_1(s)-y_1(s)| \leq
\frac{2}{m}+\frac{2}{m e_m} \mathbb{E}|h_1-y_1| \leq
\frac{2}{m}+\frac{2C}{m e_m} \sqrt{\Delta}
\end{eqnarray*}
Therefore,
\begin{eqnarray*}
\mathbb{E}|x_1(t)-h_1(t)|  \leq e_{m-1} + \frac{4 \sigma_1^2}{m} +
C \frac{\sqrt{\Delta}}{me_m} + C \int_0^t \mathbb{E}
(|x_1(s)-h_1(s)|+|x_2(s)-h_2(s)|)ds.
\end{eqnarray*}
Now, we do exact the same for $|x_2(s)-h_2(s)|$, adding the
results and setting $u(t) = |x_1(s)-h_1(s)|+|x_2(s)-h_2(s)|$ we
arrive at
\begin{eqnarray*}
\mathbb{E} u(t) \leq 2e_{m-1} + \frac{4
(\sigma_1^2+\sigma_2^2)}{m} + C \frac{\sqrt{\Delta}}{me_m} + C
\int_0^t \mathbb{E} u(s) ds
\end{eqnarray*}
Use now  Gronwall's inequality and substitute  in (19) and then
again Gronwall's inequality we arrive at
\begin{eqnarray*}
\mathbb{E}v^2(t) \leq C \sqrt{\Delta} + C
\frac{\sqrt{\Delta}}{me_m} + e_{m-1}.
\end{eqnarray*}
Choosing $m = \sqrt{\ln n^{\frac{1}{3}}}$ we deduce that
\begin{eqnarray*}
\mathbb{E}v^2(t) \leq C \frac{1}{\sqrt{\ln n}}
\end{eqnarray*}
 But
\begin{eqnarray*}
 \mathbb{E}(|x_1(t)-y_1(t)|^2 + |x_2(t)-y_2(t)|^2)  & \leq &  2
\mathbb{E}v^2(t) + 2
\mathbb{E}(|h_1(t)-y_1(t)|^2+|h_2(t)-y_2(t)|^2) \\ & \leq &  C
\frac{1}{\sqrt{\ln n}}
\end{eqnarray*}
 \qed

 \subsection{Polynomial order of convergence}
 In this section we will prove that the order of convergence is at
 least $1/2$ under further conditions on parameters and to do this we state first a proposition in which we
 show that the true solutions $x_1,x_2$ has exponential inverse moment
 bounds.

Consider the following CIR processes.
\begin{eqnarray*}
z_1(t) & = & x_1(0) + \int_0^t (k - \lambda_{11} z_1(s)) ds +
\sigma_1
\int_0^t \sqrt{z_1(s)} dW^1_s, \\
z_2(t) & = & x_2(0) + \int_0^t (l - \lambda_{21} z_2(s)) ds +
\sigma_2 \int_0^t \sqrt{z_2(s)} dW^2_s.
\end{eqnarray*}
{\bf Assumption B} Assume that there exists some strictly positive
constants $L_1(\sigma_1,k),L_2(\sigma_2,l)$
 such that
 \begin{eqnarray*}
 \mathbb{E} exp \left( \int_0^T  \frac{L_1}{z_1(s)}ds \right) < \infty,
 \quad \mathbb{E} exp \left( \int_0^T  \frac{L_2}{z_2(s)}ds
 \right)< \infty.
\end{eqnarray*}
One can see \cite{Hutz}, \cite{Kuznetsov}, \cite{Berkaoui2} for a
discussion on this assumption.

 \begin{proposition}
 Suppose that Assumption A and B hold.
Then, the following bounds are true,
\begin{eqnarray*}
 \mathbb{E} exp \left( \int_0^T  \frac{L_1}{x_1(s)}ds \right) < \infty,
 \quad \mathbb{E} exp \left( \int_0^T  \frac{L_2}{x_2(s)}ds
 \right)< \infty.
\end{eqnarray*}
\end{proposition}
\proof From the comparison theorem (see \cite{karatzas}, prop.
5.2.18) we know that $x_1(t) \geq z_1(t)$ with
\begin{eqnarray*}
z_1(t) = x_1(0) + \int_0^t (k - \lambda_{11}z_1(s))ds + \sigma_1
\int_0^t \sqrt{z_1(s)} dW_s^1.
\end{eqnarray*}
Therefore, since for $z_1$ we have exponential inverse moment
bounds  we take the result. The same holds for $x_2$. \qed

\begin{proposition}
Assume assumptions A and B. If
$\frac{L_1}{4(\sigma_1^2+\sigma_2^2)} \geq 1$ and
$\frac{L_2}{4(\sigma_1^2+\sigma_2^2)} \geq 1$
 the following rate of convergence holds, assuming
that $x_0 \in \mathbb{R}_+$,
\begin{eqnarray*}
\mathbb{E} | x_1(t) - y_1(t) |^2+|x_2(t)-y_2(t)|^2 \leq C  \Delta.
\end{eqnarray*}
\end{proposition}

\proof

Define the processes,
\begin{eqnarray*}
\gamma_1(t) &  = &\int_0^t
\frac{ds}{(\sqrt{x_1(s)}+\sqrt{h_1(s)})^2}
\\
\gamma_2(t) & = & \int_0^t
\frac{ds}{(\sqrt{x_2(s)}+\sqrt{h_2(s)})^2}
\\
 \gamma(t) &  = & \gamma_1(t) + \gamma_2(t)
\end{eqnarray*}
and then the stopping times defined by
\begin{eqnarray*}
\tau^1_l & = & \inf \{ s \in [0,T]: 4(\sigma_1^2+\sigma^2_2)
\gamma_1(s)+\frac{K}{2}s \geq l \}, \\
\tau^2_l & = & \inf \{ s \in [0,T]: 4(\sigma_1^2+\sigma^2_2)
\gamma_2(s)+\frac{K}{2}s \geq l \}, \\
\tau_l & = & \inf \{ s \in [0,T]: 4(\sigma_1^2+\sigma^2_2)
\gamma(s)+Ks \geq l \}
\end{eqnarray*}
for some fixed $K > 0$.

 Using Ito's formula on $|x_1({\tau})-h_1({\tau})|^2$ with $\tau$
a stopping time, we obtain
\begin{eqnarray*}
\mathbb{E}(x_1({\tau})-h_1({\tau}))^2 & \leq &  \int_0^{\tau}
\big(\mathbb{E} 2 \lambda_{11}|x_1(s)-h_1(s)||y_1(s)-x_1(s)| \\
& & + 2\lambda_{12} |x_1(s) - h_1(s)||x_2(s)-y_2(\hat{s})|
+\sigma_1^2 \mathbb{E}| \sqrt{y_1(s)} - \sqrt{x_1(s)}|^2 \big) ds
\\ & \leq & \int_0^{\tau} 3 \lambda_{11}
\mathbb{E}|x_1(s)-h_1(s)|^2 +2 \lambda_{11}
\mathbb{E}|h_1(s)-y_1(s)|^2 \\ & &  +    \lambda_{12}
\mathbb{E}|x_1(s)-h_1(s)|^2 + 2\lambda_{12}
\mathbb{E}|x_2(s)-h_2(s)|^2 \\  & & +2 \lambda_{12}
\mathbb{E}|h_2(s) - y_2(\hat{s})|^2+\sigma_1^2 \mathbb{E}|
\sqrt{y_1(s)} - \sqrt{x_1(s)}|^2 )  ds \\ & \leq &  C \Delta +
\int_0^{\tau} (3 \lambda_{11}+\lambda_{12}) \mathbb{E}
(|x_1(s)-h_1(s)|^2+|x_2(s)-h_2(s)|^2) \\ & & + \sigma_1^2
\mathbb{E}| \sqrt{y_1(s)} - \sqrt{x_1(s)}|^2   ds
\end{eqnarray*}

The last term can be expressed as
\begin{eqnarray*}
\int_0^{\tau} \mathbb{E}|\sqrt{y_1(s)}-\sqrt{x_1(s)}|^2 ds  & \leq
&
 \int_0^{\tau} 2\mathbb{E} |\sqrt{x_1(s)}-\sqrt{h_1(s)}|^2 + 2
\mathbb{E}|\sqrt{h_1(s)}-\sqrt{y_1(s)}|^2 ds \\ & \leq & 2
\int_0^{\tau} \mathbb{E} |\sqrt{x_1(s)}-\sqrt{h_1(s)}|^2 ds +
C\Delta.
\end{eqnarray*}
But
\begin{eqnarray*}
\int_0^{\tau} \mathbb{E} |\sqrt{x_1(s)}-\sqrt{h_1(s)}|^2 ds =
\int_0^{\tau} \mathbb{E}
\frac{|x_1(s)-h_1(s)|^2}{(\sqrt{x_1(s)}+\sqrt{h_1(s)})^2} ds.
\end{eqnarray*}

Doing exactly the same work on $|x_2(\tau)-h_2(\tau)|^2$, adding
the results and setting $v^2(\tau)  = |x_1(\tau)-h_1(\tau)|^2 +
|x_2(\tau) - h_2(\tau)|^2$ we get,
\begin{eqnarray}
\mathbb{E} v^2(\tau) & \leq & C \Delta + \int_0^{\tau}
\mathbb{E}(K v^2(s) +
\frac{2\sigma_1^2|x_1(s)-h_1(s)|^2}{(\sqrt{x_1(s)} +
\sqrt{h_1(s)})^2} +
\frac{2\sigma_2^2|x_2(s)-h_2(s)|^2}{(\sqrt{x_2(s)} +
\sqrt{h_2(s)})^2}) ds \nonumber \\ & \leq & C \Delta +
\int_0^{\tau} \mathbb{E} (Ks+4(\sigma^2_1+\sigma_2^2)
\gamma_s)^{'} v_s^2 ds
\end{eqnarray}

Now, for $\tau = \tau_l$, we use the change of variables setting
$u =4(\sigma_1^2+\sigma^2_2) \gamma(s)+Ks$ and therefore $s =
\tau_u$ obtaining,
\begin{eqnarray*}
\mathbb{E} v^2_{\tau_l} \leq C \Delta +  \int_0^l \mathbb{E}
v^2_{\tau_u} du.
\end{eqnarray*}
Using Gronwall's inequality we obtain,
\begin{eqnarray}
\mathbb{E}v^2_{\tau_l} \leq C e^l \Delta.
\end{eqnarray}
Now we rewrite  (20) as follows,
\begin{eqnarray}
\mathbb{E}v^2(\tau) & \leq &  C \Delta + \int_0^{\tau}
\mathbb{E}(\frac{K}{2}s + 4(\sigma^2_1 +
\sigma_2^2)\gamma_1(s))^{'} v^2(s) ds \nonumber  \\ & & +
\int_0^{\tau} \mathbb{E}(\frac{K}{2}s + 4(\sigma^2_1 +
\sigma_2^2)\gamma_2(s))^{'} v^2(s) ds
\end{eqnarray}
For $\tau = t \wedge \tau_l \in [0,T]$ in (22), we have under the
change of variables $u = 4(\sigma_1^2+\sigma^2_2)
\gamma_1(s)+\frac{K}{2}s$, for the first integral, and the change
of variables $u = 4(\sigma_1^2+\sigma^2_2)
\gamma_2(s)+\frac{K}{2}s$ for the second integral,
\begin{eqnarray}
\mathbb{E} v^2(t \wedge \tau_l) & \leq & C\Delta +\mathbb{E}
\int_0^{\frac{K}{2}T + 4(\sigma_1^2+\sigma^2_2)\gamma_1(T)}
v^2(\tau^1_u \wedge \tau_u) du \nonumber \\ & & +
\mathbb{E} \int_0^{\frac{K}{2}T +4(\sigma_1^2+\sigma^2_2)\gamma_2(T)} v^2(\tau^2_u \wedge \tau_u) du\nonumber \\
& \leq & C \Delta + \int_0^{\infty} \mathbb{E} \left( \mathbb{I}_{
\{ \frac{K}{2}T + 4(\sigma_1^2+\sigma^2_2) \gamma_1(T) \geq u \} }
v^2(\tau^1_u \wedge \tau_u) \right) du \nonumber \\ & & +
\int_0^{\infty} \mathbb{E} \left( \mathbb{I}_{ \{ \frac{K}{2}T +
4(\sigma_1^2+\sigma^2_2) \gamma_2(T) \geq u \} } v^2(\tau^2_u
\wedge \tau_u) \right) du.
\end{eqnarray}
Noting that
\begin{eqnarray*}
& & \int_0^{\infty} \mathbb{E} \left( \mathbb{I}_{ \{ \frac{K}{2}T
+
4(\sigma_1^2+\sigma^2_2) \gamma_1(T) \geq u \} } v^2(\tau^1_u \wedge \tau_u) \right) du \\
& \leq  & \int_0^{\frac{K}{2}T} \mathbb{E} v^2(\tau^1_u \wedge
\tau_u)du \\ & & + \int_{\frac{K}{2}T}^{\infty} \mathbb{P}
(\frac{K}{2}T + 4(\sigma_1^2+\sigma^2_2) \gamma_1(T) \geq u)
\mathbb{E} \left( v^2(\tau^1_u \wedge \tau_u) \; | \; \{
\frac{K}{2}T +
4(\sigma_1^2+\sigma^2_2) \gamma_1(T) \geq u \} \right) du,   \\
& \leq & C \Delta + \int_{0}^{\infty} \mathbb{P} (
4(\sigma_1^2+\sigma^2_2) \gamma_1(T)  \geq u) \mathbb{E}
v^2(\tau^1_u \wedge \tau_u) du
\end{eqnarray*}
and then, with exactly the same arguments for the integral
involving $\gamma_2(t)$, we arrive using (21)
\begin{eqnarray*}
\mathbb{E}v^2(t \wedge \tau_l) \leq C\Delta \left(1+  C
\int_0^{\infty} \mathbb{P} ( 4(\sigma_1^2+\sigma^2_2) \gamma_1(T)
\geq u) e^{u} du + \int_0^{\infty} \mathbb{P} (
4(\sigma_1^2+\sigma^2_2) \gamma_2(T)  \geq u) e^{u} du\right).
\end{eqnarray*}
The probability,
\begin{eqnarray*}
\mathbb{P}( 4(\sigma_1^2+\sigma^2_2) \gamma_1(T)  \geq u) \leq
\frac{1}{e^{mu}} \mathbb{E}(e^{4m(\sigma_1^2+\sigma^2_2)
\gamma_1(T) } ),
\end{eqnarray*}
and the same holds for the probability involving $\gamma_2(t)$.
Choose $m_i = \frac{L_i}{4(\sigma_1^2+\sigma_2^2)}$, for $i=1,2$
and use Proposition 3 to deduce that
\begin{eqnarray*}
\mathbb{E} v^2(t \wedge \tau_l) \leq C \Delta.
\end{eqnarray*}
Using Fatou's lemma for $l \to \infty$ we take the result.
 \qed

\section{A second explicit numerical scheme}
 We will propose a
different numerical scheme below,
\begin{eqnarray*}
y_1(t_{k+1}) & = & \left(
\frac{\sigma_1}{2}(W^1_{t_{k+1}}-W^1_{t_{k}}) +
\sqrt{y_1(t_k)(1-\lambda_{11} \Delta)+ \Delta \lambda_{12}
y_2(t_k) + \Delta (k-
\frac{\sigma_1^2}{4})} \right)^2 , \\
y_2(t_{k+1}) & = & \left(
\frac{\sigma_2}{2}(W^2_{t_{k+1}}-W^2_{t_{k}}) +
\sqrt{y_2(t_k)(1-\lambda_{21}\Delta)+ \Delta \lambda_{22} y_1(t_k)
+ \Delta (l- \frac{\sigma_2^2}{4})} \right)^2.
\end{eqnarray*}
Knowing $y_1(t_0) = x_1(0), y_2(t_0) = x_2(0)$ we obtain
explicitly and parallel the $y_1(t_1), y_2(t_2)$ and so on.

We work with the following stochastic processes,
\begin{eqnarray*}
y_1(t) & = & \left( \frac{\sigma_1}{2}(W^1_{t}-W^1_{t_{k}}) +
\sqrt{y_1(t_k)(1-\lambda_{11} \Delta)+ \Delta \lambda_{12}
y_2(t_k) + \Delta (k-
\frac{\sigma_1^2}{4})} \right)^2 = (z_1(t))^2, \\
y_2(t) & = & \left( \frac{\sigma_2}{2}(W^2_{t}-W^2_{t_{k}}) +
\sqrt{y_2(t_k)(1-\lambda_{21}\Delta)+ \Delta \lambda_{22} y_1(t_k)
+ \Delta (l- \frac{\sigma_2^2}{4})} \right)^2 = (z_2(t))^2,
\end{eqnarray*}
and in fact with the stochastic differentials obtained by the use
of Ito's formula, for $t \in (t_k,t_{k+1}]$,
\begin{eqnarray*}
y_1(t) & = & y_1(t_k)(1-\lambda_{11} \Delta)+ \Delta \lambda_{12}
y_2(t_k) + \Delta (k- \frac{\sigma_1^2}{4}) + \int_{t_k}^t
\frac{\sigma_1^2}{4} ds  \\ & & + \sigma_1 \int_{t_k}^t
sgn(z_1(s)) \sqrt{y_1(s)}dW^1_s, \end{eqnarray*} \begin{eqnarray*}
y_2(t) & = & y_2(t_k)(1-\lambda_{21} \Delta)+ \Delta \lambda_{22}
y_1(t_k) + \Delta (l- \frac{\sigma_2^2}{4}) + \int_{t_k}^t
\frac{\sigma_2^2}{4} ds \\ & & + \sigma_2 \int_{t_k}^t sgn(z_2(s))
\sqrt{y_2(s)}dW^2_s.
\end{eqnarray*}

The compact forms are, for $t \in (t_k,t_{k+1}]$,
\begin{eqnarray*}
y_1(t)   =  & & x_1(0)+\int_0^t(k-\lambda_{11} y_1(\hat{s}) +
\lambda_{12} y_2(\hat{s})) ds + \int_t^{t_{k+1}} (k -
\frac{\sigma_1^2}{4} -\lambda_{11} y_1(t_k) + \lambda_{12}
y_2(t_k))ds \\ &  + & \sigma_1 \int_0^t sgn(z_1(s)) \sqrt{y_1(s)}
dW^1_s,
\\
y_2(t)  = & & x_2(0)+\int_0^t(l-\lambda_{21} y_2(\hat{s}) +
\lambda_{22} y_1(\hat{s})) ds + \int_t^{t_{k+1}} (l -
\frac{\sigma_2^2}{4} -\lambda_{21} y_2(t_k) + \lambda_{22}
y_1(t_k))ds \\ &  + & \sigma_2 \int_0^t sgn(z_2(s)) \sqrt{y_2(s)}
dW^2_s.
\end{eqnarray*}

Finally, we will use the following auxiliary processes,
\begin{eqnarray*}
h_1(t) &  = & x_1(0)+\int_0^t(k-\lambda_{11} y_1(\hat{s}) +
\lambda_{12} y_2(\hat{s})) ds   + \sigma_1 \int_0^t sgn(z_1(s))
\sqrt{y_1(s)} dW^1_s,
\\
h_2(t) & = & x_2(0)+\int_0^t(l-\lambda_{21} y_2(\hat{s}) +
\lambda_{22} y_1(\hat{s})) ds + \sigma_2 \int_0^t sgn(z_2(s))
\sqrt{y_2(s)} dW^2_s.
\end{eqnarray*}

{\bf Assumption C} Assume that $d_1 \geq 1$, $d_2 \geq 1$, $\Delta
\leq \frac{1}{\max \{\lambda_{11}, \lambda_{21} \}}$ and $x_0 \in
\mathbb{R}_+$.

\begin{lemma}
Under Assumption C we have
\begin{eqnarray*}
\mathbb{E}(\sup_{0 \leq t \leq T} (y_1(t)^2 + y_2(t)^2) )< C
\end{eqnarray*}
\end{lemma}

\proof Here, again, we easily see that
\begin{eqnarray*}
0 \leq y_1(t) \leq v_1(t) = x_1(0) & + & Tk + \Delta \lambda_{12}
x_2(0)+\int_0^t \lambda_{12} (y_2(\hat{s})+y_2(t_k)) ds \\ &  + &
\sigma_1 \int_0^t \sqrt{y_1(s)} dW^1_s,\end{eqnarray*}
\begin{eqnarray*} 0 \leq y_2(t) \leq v_2(t) = x_2(0) & + & Tl + \Delta
\lambda_{22} x_1(0)+\int_0^t \lambda_{22} (y_1(\hat{s})+y_1(t_k))
ds \\ & + & \sigma_2\int_0^t \sqrt{y_2(s)} dW^2_s,
\end{eqnarray*}
Continuing as before we get the result. \qed

\begin{lemma} Under Assumption C we have the following estimates,
for $i=1,2$ and $t \in [t_k,t_{k+1}]$,
\begin{eqnarray*}
\mathbb{E}|h_i(t) - y_i(t)|^2 & \leq & C \Delta \\
\mathbb{E} |y_i(t) - y_i(t_k)|^2 & \leq & C \Delta, \\
\mathbb{E} |h_i(t) - y_i(t_k)|^2 & \leq & C \Delta,\\
\mathbb{E} |h_i(t)|^2 & \leq & C \Delta.
\end{eqnarray*}
\end{lemma}

\proof Using the moment bounds of Lemma 8 we easily get the
result. \qed

\begin{lemma}
Under Assumption B, we have the following estimates,
\begin{eqnarray*}
\mathbb{E} y_1(t) \left(sgn(z_1(t)) - 1\right)^2 \leq C
\Delta^{\frac{1}{2} - \vep}, \quad \mathbb{E} y_2(t)
\left(sgn(z_2(t)) - 1\right)^2 \leq C \Delta^{\frac{1}{2} - \vep}
\end{eqnarray*}
for any $\vep > 0$.
\end{lemma}

\proof We begin with, when $t \in [t_k,t_{k+1}]$,
\begin{eqnarray*}
\mathbb{E} y_1(t) (sgn(z_1(t)) - 1)^2 & = & 4\mathbb{E}y_1(t)
\mathbb{I}_{ \{
 z_1(t) \leq 0 \} } \leq 4 \mathbb{E} |y_1(t)-y_1(t_k)| + 4 \mathbb{E} y_1(t_k) \mathbb{I}_{
\{
 z_1(t) \leq 0 \} }   \\ & \leq & C \Delta +4 \mathbb{E} y_1(t_k) \mathbb{I}_{
\{
 z_1(t) \leq 0 \} } \mathbb{I}_{ \{ y_1(t_k) \leq \Delta^{1/2-\vep} \} } +
 4 \mathbb{E} y_1(t_k) \mathbb{I}_{
\{
 z_1(t) \leq 0 \} } \mathbb{I}_{ \{ y_1(t_k) > \Delta^{1/2-\vep} \} } \\
 & \leq & C \Delta^{1/2-\vep} + 4 \mathbb{E} y_1(t_k) \mathbb{I}_{
\left\{
 \{z_1(t) \leq 0 \} \cap \{  y_1(t_k) > \Delta^{1/2-\vep} \} \right\} }
 \end{eqnarray*}
 We have used Lemma 9 to obtain the second inequality, estimating  the term $\mathbb{E} |y_1(t)-y_1(t_k)|$.
But
\begin{eqnarray*}
  & & \{ z_1(t) \leq 0\} \cap \{ y_1(t_k) > \Delta^{1/2-\vep } \}
  \\
& = & \left\{ W^1_t - W^1_{t_k} \leq -\frac{2 }{\sigma_1}
\sqrt{y_1(t_k)(1-\lambda_{11}\Delta)
 +\Delta \lambda_{12}y_2(t_k)+ \Delta (k - \frac{\sigma_1^2}{4}) } \right\} \cap \{ y_1(t_k) > \Delta^{1/2-\vep } \} \\
& \subseteq & \left\{ W^1_t - W^1_{t_k} \leq
-\frac{2\sqrt{1-\lambda_{11} \Delta}}{\sigma_1}
 \sqrt{\Delta^{1/2-\vep}} \right\} .
\end{eqnarray*}
 Since the increment $W^1_t-W^1_{t_k}$ is normally distributed
with mean zero and variance $t-t_k$ we  have that
\begin{eqnarray*}
\mathbb{P}\left( \{ z_1(t) \leq 0 \} \cap \{ y_1(t_k) >
\Delta^{1/2-\vep } \} \right) \leq C\frac{\sqrt{t-t_k}}{\sqrt{2
\pi (t-t_k)}} \int_{C\frac{2
\sqrt{\Delta^{1/2-\vep}}}{\sqrt{t-t_k}}}^{\infty} e^{-y^2/2} dy
\leq  \frac{C\Delta^{\vep}}{e^{C/\Delta^{\vep}}}.
\end{eqnarray*}
We have used the inequality of problem 9.22, p.112 of
\cite{karatzas} to obtain the last inequality. Now we have, using
the moment bounds for the numerical solution,
\begin{eqnarray*}
\mathbb{E} y_1(t_k) \mathbb{I}_{ \left\{
 \{z_1(t) \leq 0 \} \cap \{  y_1(t_k) > \Delta^{1/2-\vep} \} \right\}
 }  \leq C \mathbb{P}\left( \{ z_1(t) \leq 0 \} \cap \{ y_1(t_k) >
\Delta^{1/2-\vep } \} \right)
\end{eqnarray*}
Noting that  $\frac{\Delta^{\vep}}{e^{1/\Delta^{\vep}}} \to 0$
faster than any power of $\Delta$  we have that
\begin{eqnarray*}
\mathbb{E} y_1(t) (sgn(z_1(t)) - 1)^2 \leq C \Delta^{\frac{1}{2} -
\vep}.
\end{eqnarray*}
The same holds for $y_2(t)$.
 \qed

Because $h_1,h_2$ are essential the same as in the previous
section,  we can use the same arguments as in Theorem 3 and
Proposition 4 together with Lemma 10 to get the following results.
\begin{theorem}
If Assumption C holds then
\begin{eqnarray*}
 \mathbb{E} (|x_1(t) - y_1(t)|^2+|x_2(t)-y_2(t)|^2)  \leq C \frac{1}{\sqrt{\ln n}}
\end{eqnarray*}
for any $t \in [0,T]$.
\end{theorem}

\begin{proposition}
Suppose that Assumptions B and C hold. Then, if
$\frac{L_1}{4(\sigma_1^2+\sigma_2^2)} \geq 1$ and
$\frac{L_2}{4(\sigma_1^2+\sigma_2^2)} \geq 1$,
 the following rate of convergence holds,
\begin{eqnarray*}
\mathbb{E} | x_1(t) - y_1(t) |^2+|x_2(t)-y_2(t)|^2 \leq C
\Delta^{1/2-\vep}.
\end{eqnarray*}
for every $\vep > 0$. That is the order of convergence is at least
$1/4-\vep$.
\end{proposition}

{\bf Conclusion} We have proposed two explicit and positivity
preserving numerical schemes for the two factor CIR model. The
first one is based on the exact simulation of the CIR process for
a specific set of parameters. The advantage of the second method
is that one need less calculations in each step comparing with the
first method. However, extended numerical experiments has to be
done to compare them. Let us mention that both the results hold
for the case of one equation choosing for example $\lambda_{12} =
0$. Finally, the above results can be easily   extended for the
multi-factor case.

In \cite{Wilkie} one can find a different use of the above model.
If one considers a more complicated model than the above, for
example,
\begin{eqnarray*}
x_1(t) & = & x_1(0) + \int_0^t (k - \lambda_{11} x_1(s) +
\lambda_{12}x_2(s))ds +
\int_0^t \sigma_1 \sqrt{x_1(s)x_2(s)} dW^1_s,\\
x_2(t) & = & x_2(0) + \int_0^t (l - \lambda_{21} x_2(s) +
\lambda_{22} x_1(s))ds + \int_0^t \sigma_2 \sqrt{x_1(s)x_2(s)}
dW^2_s
\end{eqnarray*}
then it is not obvious how our first numerical scheme based on
exact simulation can be applied here. Considering the second
method one can propose the following numerical scheme, {\small
\begin{eqnarray*} y_1(t_{k+1}) & = & \left( \frac{\sigma_1
\sqrt{y_2(t_k)}}{2}(W^1_{t_{k+1}}-W^1_{t_{k}}) +
\sqrt{y_1(t_k)(1-\lambda_{11} \Delta)+ \Delta \lambda_{12}
y_2(t_k) + \Delta (k-
\frac{\sigma_1^2 y_2(t_k)}{4})} \right)^2 , \\
y_2(t_{k+1}) & = & \left(
\frac{\sigma_2\sqrt{y_1(t_k)}}{2}(W^2_{t_{k+1}}-W^2_{t_{k}}) +
\sqrt{y_2(t_k)(1-\lambda_{21}\Delta)+ \Delta \lambda_{22} y_1(t_k)
+ \Delta (l- \frac{\sigma_2^2y_1(t_k)}{4})} \right)^2.
\end{eqnarray*} }
With the same analysis and with a minor modification on the
hypotheses, one can prove that this scheme converges strongly to
the true solution but without some rate, i.e. a similar result as
Theorem 3.

As a minimal computer experiment we give below the difference
between the numerical scheme (2) for $a=1$ and the scheme proposed
in \cite{Alfonsi} just to see that these methods are close. More
complicated computer experiments has to be done in order to detect
the actual order of convergence and other advantages or
disadvantages of this method compared with that of \cite{Alfonsi}.

\begin{figure}[ht]
  \caption{    $x_0 = 4$, $\Delta=10^{-4}$, $k=2$, $l=1$, $s=1$ $T=1$.}
  \centering
    \includegraphics[width=0.3\textwidth]{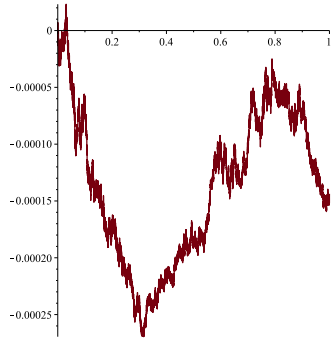}
\end{figure}


\section*{Acknowledgement}

I would like to thank Prof. Arnulf Jentzen for an instructive
discussion during the preparation of this paper.

\end{document}